\documentclass[11pt]{article}

\usepackage{amsmath}
\usepackage{amssymb}
\usepackage{graphics}
\usepackage{rotating}
\usepackage{cite}
\usepackage{color}

\def\0{\mbox{\tiny $0$}}
\def\1{\mbox{\tiny $1$}}
\def\2{\mbox{\tiny $2$}}
\def\3{\mbox{\tiny $3$}}
\def\4{\mbox{\tiny $4$}}
\def\5{\mbox{\tiny $5$}}
\def\6{\mbox{\tiny $6$}}
\def\7{\mbox{\tiny $7$}}
\def\8{\mbox{\tiny $8$}}
\def\9{\mbox{\tiny $9$}}
\def\i{\mbox{\tiny $i$}}
\def\j{\mbox{\tiny $j$}}
\def\k{\mbox{\tiny $k$}}
\def\m{\mbox{\tiny $m$}}
\def\n{\mbox{\tiny $n$}}
\def\p{\mbox{\tiny $p$}}

\def\r{\mbox{\tiny $r$}}
\def\s{\mbox{\tiny $s$}}
\def\({\mbox{\tiny $($}}
\def\){\mbox{\tiny $)$}}

\def\xt{\mbox{\tiny $\times$}}

\def\mi{\mbox{\tiny $-$}}

\def\={\mbox{\tiny $=$}}
\def\Re{\mbox{\tiny $\mathbb{R}$}}
\def\Co{\mbox{\tiny $\mathbb{C}$}}
\def\Ha{\mbox{\tiny $\mathbb{H}$}}
\def\Xf{\mbox{\tiny $\mathbb{X}$}}
\def\mup{\mbox{\tiny $\mu$}}
\def\nup{\mbox{\tiny $\nu$}}
\begin{document}
%

\title{{\huge {\bf Real linear quaternionic operators}}}
\author{
{\Large \sf Stefano De Leo$^{(a)}$ {\sf and} Gisele Ducati$^{(b)}$}\\ \\
{\small $^{(a)}${\em Departamento de Matem\'atica Aplicada, UNICAMP}}\\
{\small  CP 6065, 13081-970  Campinas (SP) Brasil}\\
{\small  \tt deleo/ducati@ime.unicamp.br}\\
{\small  $^{(b)}${\em Departamento de Matem\'atica, UFPR}}\\
{\small  CP 19081, 81531-970  Curitiba (PR) Brasil}\\
{\small  \tt ducati@mat.ufpr.br} }

\maketitle



\abstract{In a recent paper [J.Math.Phys. {\bf 42}, 2236--2265
(2001)], we discussed differential operators within a quaternionic
formulation of quantum mechanics. In particular, we proposed a
practical method to solve quaternionic and complex linear second
order differential equations with constant coefficients. In this
paper, we extend our discussion to real linear quaternionic
differential equations. The method of resolution is based on the
Jordan canonical form of quaternionic matrices associoated to real
linear differential
operators.\\
~\\
~\\
{\bf MSC.} 15A18 -- 15A21 -- 15A30 -- 15A90 -- 47D25 -- 47E05.
}













\maketitle


\section*{I. Introduction}

Before going into the details of the discussion of $\mathbb{R}$
linear quaternionic differential operators, we briefly recall the
technique used to solve $\mathbb{H}$ and $\mathbb{C}$ linear
differential equations with constant quaternionic coefficients and
show the difficulties in extending the method of resolution to the
$\mathbb{R}$-linear case. We do not attempt a formal discussion of
quaternionic theory of differential equations. Instead, we take an
operational and intuitive approach. For the convenience of the
reader and to make our exposition as self-contained as possible,
we follow the mathematical notation and terminology used in our
previous paper~\cite{DeDu01}. In particular, the operators
\[
L_{\mup} = (1, L_{\i}, L_{\j}, L_{\k})~~~\mbox{and}~~~
R_{\mup} = (1,
R_{\i}, R_{\j}, R_{\k})~,~~~{\mbox{\small $\mu = 0,1,2,3$}}~,
\]
will denote  the left and right action  on quaternionic
functions of real variable, $\varphi(x)$,
of the
imaginary units $i, j$ and $k$. To shorten notation, we
shall use the upper-script $\mathbb{X} = \{ \mathbb{R},
\mathbb{C}, \mathbb{H} \}$ to indicate the $\mathbb{R}$, $\mathbb{C}$,
and $\mathbb{H}$ linearity (from the right) of $n$-order
quaternionic differential operators, $\mathcal{D}_{\n}^{\Xf}$,
and $n$-dimensional quaternionic matrices, $M_{\n}^{\Xf}$.

The general solution of $n$-order $\mathbb{C}$ linear homogeneous ordinary
differential  equations
\begin{equation}
\label{de}
\mathcal{D}_{\n }^{\Co} \, \varphi(x) = 0~,
\end{equation}
where
\[
\mathcal{D}_{\n}^{\Co}  =
\frac{d^{^{\, \n}}}{dx^{^{\n}}} \, - \sum_{\p \= \0}^{\n \mi \1}
a^{\( \p \)}_{\Co}
\, \frac{d^{^{\, \p}}}{dx^{\p}} =
\frac{d^{^{\, \n}}}{dx^{^{\n}}} \, - \sum_{\p \= \0}^{\n \mi \1}
\left[
\sum_{\mup  \= \0}^{\3}
\sum_{\nup  \= \0}^{\1}
a^{\( \p \)}_{\mup \nup}  L_{\mup} R_{\nup}
\right]
\, \frac{d^{^{\, \p}}}{dx^{\p}}~,~~~
a^{\( \p \)}_{\mup \nup} \in \mathbb{R}~,
\]
has the form
\begin{equation}
\label{gsol}
\varphi(x) = \sum_{\m \= \1}^{\2 \n} \varphi_{\m} (x) \, c_{\m}~,
\end{equation}
where $\left\{ \varphi_{\1} (x) \, , \, \dots \, , \,
\varphi_{\2 \n} (x) \right\}$ represent  $2n$  quaternionic
particular solutions linearly independent
over $\mathbb{C}(1,i)$ and $c_{\m}$ are
complex constants determined by the initial values of the function
$\varphi(x)$ and its derivatives
\begin{equation}
\label{icon}
\varphi(x_{\0}) = \varphi_{\0}~,~~
\frac{d \varphi}{d x} (x_{\0})
= \varphi_{\1}~,~\dots~,~
\frac{d^{^{\n \mi \1}} \varphi}{d x^{\n \mi \1}}
(x_{\0})
= \varphi_{\n \mi \1 }~~\in~\mathbb{H}~.
\end{equation}
The solution of
Eq.~(\ref{de}) is given by
\begin{equation}
\label{sol}
\varphi(x)  =  \sum_{\p \= \1}^{\n} \left\{
\exp \left[
M_{\n }^{\Co} \, \left( x - x_{\0} \right) \right]
\right\}_{\1 \p} \, \varphi_{\p \mi \1}
 =  \sum_{\p \= \1}^{\n} \left\{ S_{\n }^{\Co} \,
\exp \left[
J_{\n }^{\Co} \, \left( x - x_{\0} \right) \right]
(S_{\n }^{\Co})^{\mi \1}
\right\}_{\1 \p} \, \varphi_{\p \mi \1}~,
\end{equation}
where
\begin{equation}
\label{mc}
M_{\n }^{\Co} =
\left( \begin{array}{lllllll}
 \, 0 & 1 & 0 & \cdot & ~\cdot & ~\cdot & \, ~0\\
\,  0 & 0 & 1 & \cdot & ~\cdot & ~\cdot & \, ~0\\
\,  \cdot & \cdot & \cdot &  \cdot & ~\cdot & ~\cdot & \, ~\cdot \\
\,  \cdot & \cdot & \cdot &  \cdot & ~\cdot & ~\cdot & \, ~\cdot \\
\,  \cdot & \cdot & \cdot &  \cdot & ~\cdot & ~\cdot & \, ~\cdot \\
\, 0 & 0 & 0 & \cdot & ~\cdot & ~\cdot & \, ~1\\
\, a_{\Co}^{\( \0 \)} &
a_{\Co}^{\( \1 \)} &
a_{\Co}^{\( \2 \)} &
\cdot & ~\cdot & ~\cdot &  \, ~a_{\Co}^{\( \n \mi \1 \)}
\end{array} \right)
\end{equation}
and $J_{\n }^{\Co}$ represents the  Jordan canonical form of
$M_{\n }^{\Co}$. Note that
$J_{\n }^{\Co}= D_{\n }^{\Co} +
N_{\n }^{\Co}$, where $D_{\n }^{\Co}$
is diagonal and $N_{\n }^{\Co}$ is nilpotent. The Jordan matrix
$J_{\n }^{\Co}$
can be determined by solving the right complex eigenvalue
problem~\cite{Zha97,Bak99,DeSco00}
for $M_{\n }^{\Co}$.

It is worth pointing out that
from Eq.~(\ref{sol}) we can recover the form of the particular solutions
$\varphi_{\m}(x)$. For example, in the case
of a null nilpotent matrix, the general solution (\ref{gsol})
can be rewritten in terms of  left acting
quaternionic coefficients ($u_{\p}$ and $v_{\p}$), complex
exponentials ($\exp \, [ \, z_{\p} \, x \, ]$ and
$\exp \, [ \, w_{\p} \, x \, ]$), and right acting complex constants
($c_{\p}$ and $\tilde{c}_{\p}$)  determined by  the initial conditions.
Explicitly, we find
\begin{equation}
\label{gsol2}
\varphi(x) = \sum_{\p \= \1}^{\n} \left\{ u_{\p} \,
\exp \, [ \, z_{\p} \, x \, ] \, c_{\p} +
v_{\p} \,
\exp \, [ \, w_{\p} \, x \, ] \, \tilde{c}_{\p}
\right\}~,
\end{equation}
where the  complex coefficients
$\left\{ z_{\1} \, , \, w_{\1} \, , \, \dots \, , \,
z_{\n} \, , \, w_{\n} \right\}$ represent the right
eigenvalues of  the quaternionic matrix
$M_{\n \, }^{\Co}$. In the case of equal complex eigenvalues
and equal quaternionic coefficients the remaining particular
solutions are determined by using
the   nilpotent matrix $N_{\n}^{\Co}$.

For $\mathbb{H}$ linear quaternionic matrices,
$M_{\n}^{\Ha}$, it is possible to show that $v_{\p}=u_{\p} \, j$ and
$w_{\p}=\bar{z}_{\p}$. Consequently,  for
$\mathbb{H}$ linear quaternionic differential operators,
the general solution (\ref{gsol2}) reduces to
\begin{equation}
\label{gsol3}
\varphi(x)  =  \sum_{\p \= \1}^{\n} u_{\p} \,
\exp \, [ \, z_{\p} \, x \, ] \, (c_{\p} +
j \tilde{c}_{\p} )
  =  \sum_{\p \= \1}^{\n}
\exp \, [ \, q_{\p} \, x \, ] \, h_{\p}~,
\end{equation}
where $q_{\p} = u_{\p} z_{\p} u_{\p}^{\mi \1}$. The initial
conditions (\ref{icon}) shall fix
the $n$  quaternionic
constants $h_{\p}$.

Due to the
$\mathbb{R}$ linearity, the general solution of $n$-order
homogeneous ordinary differential  equations with quaternionic constant
coefficients  which appear on the left and on the right has the form
\begin{equation}
\label{gsolr}
\varphi(x) = \sum_{\s \= \1}^{\4 \n} \varphi_{\s} (x) \, r_{\s}~,
\end{equation}
where $\left\{ \varphi_{\1} (x) \, , \, \dots \, , \,
\varphi_{\4 \n} (x) \right\}$ represent  $4n$  quaternionic
particular solutions, linearly independent
over $\mathbb{R}$, and $r_{\s}$ are real constants fixed by Eqs.~(\ref{icon}).
The question still unanswered is how to determine the particular solutions
$\varphi_{\s} (x)$.  The natural choice of left acting quaternionic
coefficients and real exponentials
\[
\varphi_{\s}(x) =
u_{\s} \exp \, [ \, \lambda_{\s} \, x \, ]~,
\]
does not represent a satisfactory answer. In fact, such particular solutions
are at most valid for the  real part of the eigenvalue spectrum of the matrix
$M_{\n}^{\Re}$.


\section*{II. Real linear quaternionic differential equations}

Let us consider the second
order $\mathbb{R}$ linear quaternionic
differential  equation
\begin{equation}
\label{de2}
\mathcal{D}_{\2 }^{\Re} \, \varphi(x) = 0~,
\end{equation}
where
\[
\mathcal{D}_{\2}^{\Re}   =
\frac{d^{^{\, \2}}}{dx^{^{\2}}} \, -
a^{\( \1 \)}_{\Re}
\, \frac{d}{dx} - a^{\( \0 \)}_{\Re}   =
\frac{d^{^{\, \2}}}{dx^{^{\2}}} \, -
\sum_{\mup , \nup  \= \0}^{\3}
a^{\( \1 \)}_{\mup \nup}  L_{\mup} R_{\nup}
\, \frac{d}{dx} -
\sum_{\mup , \nup  \= \0}^{\3}
a^{\( \0 \)}_{\mup \nup}  L_{\mup} R_{\nup}~.
\]
By introducing the $\mathbb{R}$ linear quaternionic matrix
\[
M_{\2}^{\Re} =
\left( \begin{array}{cc} 0 & ~1\\  & \\a^{\( \0 \)}_{\Re} &
~a^{\( \1 \)}_{\Re}
\end{array} \right)
\]
and the quaternionic column vector
\[
\Phi = \left[ \begin{array}{c} \varphi \\ \\
{\frac{d}{dx}} \, \varphi\end{array} \right]~,
\]
we can rewrite Eq.~(\ref{de2}) in matrix form
\begin{equation}
\label{matrix}
\mbox{$\frac{d}{dx}$} \, \Phi (x) = M_{\2}^{\Re} \, \Phi (x)~.
\end{equation}
The real matrix counterpart of the quaternionic
operator $M_{\2}^{\Re}$, from now on denoted by  $M_{\8}[\mathbb{R}]$, has
an 8-dimensional eigenvalue spectrum characterized by
real numbers and/or complex conjugate pairs
(the translation tables are given in the appendix).
Let $J_{\8}[\mathbb{C}]$ be  the complex Jordan form of
$M_{\8}[\mathbb{R}]$,
\begin{equation}
J_{\8}[\mathbb{C}] = R_{\2 \m} \oplus Z_{\n} \oplus
                                    \bar{Z}_{\n}~,
~~~\mbox{\small $m + n = 4$}~,
\end{equation}
where $R_{\2 \m}$ and  $Z_{\n}$ represent the matrix Jordan blocks
containing, respectively,  the real and complex eigenvalues of
$M_{\8}[\mathbb{R}]$. By using an appropriate similarity matrix
\begin{equation}
S_{\8}[\mathbb{C}] = \left(
                    \begin{array}{ccc}
 S^{\( \1 \)}_{\2 \m \mbox{\xt} \2 \m }[\mathbb{R}] &
 ~S^{\( \4 \)}_{\2 \m \mbox{\xt} \n }[\mathbb{C}]~ &
\bar{S}^{\( \4 \)}_{\2 \m \mbox{\xt} \n }[\mathbb{C}]  \\
  S^{\( \2 \)}_{\n \mbox{\xt} \2 \m }~\,[\mathbb{R}] &
 ~S^{\( \5 \)}_{\n \mbox{\xt} \n }~\,[\mathbb{C}] ~ &
\bar{S}^{\( \5 \)}_{\n \mbox{\xt} \n }~\,[\mathbb{C}]  \\
  S^{\( \3 \)}_{\n \mbox{\xt} \2 \m}~\,[\mathbb{R}] &
 ~S^{\( \6 \)}_{\n \mbox{\xt} \n }~\,[\mathbb{C}] ~ &
\bar{S}^{\( \6 \)}_{\n \mbox{\xt} \n }~\,[\mathbb{C}]
                  \end{array}    \right)~,
\end{equation}
we can rewrite  $M_{\8}[\mathbb{R}]$ as product of three complex matrices,
that is  $S_{\8}[\mathbb{C}] J_{\8}[\mathbb{C}]
S_{\8}^{\mi \1}[\mathbb{C}]$. The problem is now represented by the
impossibility to translate the single terms of the previous matrix product
by $\mathbb{R}$ linear quaternionic matrices. To overcome this difficulty we
introduce  the complex  matrix
\begin{equation}
W_{\8}[\mathbb{C}] =
\boldsymbol{1}_{\2 \m} \oplus \left[ \,  \mbox{$\frac{1}{\sqrt{2}}$} \,
\left(       \begin{array}{cr}
                      1    & i\\
                      1 & $-$ i
\end{array}   \right) \otimes \boldsymbol{1}_{\n} \, \right]~.
\end{equation}
An easy algebraic calculation shows that
\[
M_{\8}[\mathbb{R}] =
T_{\8}[\mathbb{R}] \, J_{\8}[\mathbb{R}] \,
T_{\8}^{\mi\1}[\mathbb{R}]
\]
where
\begin{equation}
T_{\8}[\mathbb{R}] = S_{\8}[\mathbb{C}] \, W_{\8}[\mathbb{C}] =  \left(
                    \begin{array}{ccc}
 S^{\( \1 \)}_{\2 \m \mbox{\xt} \2 \m }[\mathbb{R}] &
 ~~\sqrt{2} \,
\mbox{Re}\left\{ S^{\( \4 \)}_{\2 \m \mbox{\xt} \n }[\mathbb{C}] \right\}~ &
- \sqrt{2} \,
\mbox{Im}\left\{ S^{\( \4 \)}_{\2 \m \mbox{\xt} \n }[\mathbb{C}] \right\} \\
  S^{\( \2 \)}_{\n \mbox{\xt} \2 \m }~\,[\mathbb{R}] &
 ~~
\sqrt{2} \, \mbox{Re}\left\{ S^{\( \5 \)}_{\n \mbox{\xt} \n }~\,[\mathbb{C}]
\right\}~ &
- \sqrt{2} \,
\mbox{Im}\left\{ S^{\( \5 \)}_{\n \mbox{\xt} \n }~\,[\mathbb{C}] \right\} \\
  S^{\( \3 \)}_{\n \mbox{\xt} \2 \m}~\,[\mathbb{R}] &
 ~~\sqrt{2} \, \mbox{Re}\left\{ S^{\( \6 \)}_{\n \mbox{\xt} \n }~\,
[\mathbb{C}] \right\}~ &
- \sqrt{2} \,
\mbox{Im}\left\{ S^{\( \6 \)}_{\n \mbox{\xt} \n }~\,[\mathbb{C}] \right\}
                  \end{array}    \right)~,
\end{equation}
and
\begin{equation}
J_{\8}[\mathbb{R}] =
 W_{\8}^{-\1}[\mathbb{C}] \,
J_{\8}[\mathbb{C}] \, W_{\8}[\mathbb{C}]
= R_{\2 \m} \oplus \left(
                    \begin{array}{cr}
                  \mbox{Re}[Z_{\n}] & \, - \mbox{Im}[Z_{\n}]\\
                  \mbox{Im}[Z_{\n}] & \mbox{Re}[Z_{\n}]
                    \end{array}    \right)~.
\end{equation}
The real Jordan canonical form $J_{\8}[\mathbb{R}]$ can be decomposed
into the sum of  three commuting
real matrices, that is the diagonal matrix
\begin{eqnarray*}
D_{\8}[\mathbb{R}] & = &
\mbox{Diag} \, \left\{ \lambda_{\1} \, , \, ... \, , \, \lambda_{\m} \, , \,
\tilde{\lambda}_{\1} \, , \, ... \, , \, \tilde{\lambda}_{\m} \, , \,
\mbox{Re} \left[ z_{\1} \right] \, , \, ... , \, \mbox{Re} \left[ z_{\n} \right]
\, , \,
\mbox{Re} \left[ z_{\1} \right] \, , \, ... \, , \,
\mbox{Re} \left[ z_{\n} \right]
\right\}~,
\end{eqnarray*}
the anti-symmetric matrix
\begin{eqnarray*}
A_{\8}[\mathbb{R}] & = & \boldsymbol{0}_{\2m} \oplus \left[ \,
\left( \begin{array}{cr}
                                               0 & \, - 1\\
                                               1 & 0
                                    \end{array} \right) \otimes
           \mbox{Diag} \, \left\{
\mbox{Im} \left[ z_{\1} \right] \, , \, ... \, , \,
\mbox{Im} \left[ z_{\n} \right]
\right\} \, \right]~,
\end{eqnarray*}
and the nilpotent matrix $N_{\8}[\mathbb{R}]$ (a
lower triangular matrix whose only nonzero elements
are ones which appear in some of the sub-diagonal positions).
The real Jordan form
$J_{\8}[\mathbb{R}] =D_{\8}[\mathbb{R}] + A_{\8}[\mathbb{R}] +
N_{\8}[\mathbb{R}] $ and the similarity transformation
$T_{\8}[\mathbb{R}]$ can be now translated in their real linear
quaternionic counterparts
$J_{\2}^{\Re} =D_{\2}^{\Re} + A_{\2}^{\Re} + N_{\2}^{\Re}$
and $T_{\2}^{\Re}$.
The matrix solution of Eq.(\ref{matrix}) then reads
\begin{eqnarray}
\label{matsol}
\Phi(x) & = & \exp[M_{\2}^{\Re}\,  (x - x_{\0})] \Phi(x_{\0}) \nonumber \\
& = &
T_{\2}^{\Re} \, \exp[J_{\2}^{\Re} \, (x- x_{\0})] \,
( T_{\2}^{\Re})^{\mi \1} \, \Phi(x_{\0})
\nonumber \\
& = &
T_{\Re} \, \exp[D_{\Re} \, (x- x_{\0})] \,  \exp[A_{\Re} \, (x- x_{\0})] \,
\exp[N_{\Re} \, (x- x_{\0})] \,
(T_{\2}^{\Re})^{\mi \1} \, \Phi(x_{\0})~.
\end{eqnarray}
In the case of  a null nilpotent matrix (explicit examples are given in
subsection~A), it can be verified by direct calculations that
the general solution of real linear quaternionic differential equations
with constant coefficients is
\begin{eqnarray}
\label{genr}
\varphi(x) & = & \sum_{\r \= \1}^{\m} \left\{
u_{\r} \, \exp[ \lambda_{\r} \, x ] \,
\alpha_{\r} +
\tilde{u}_{\r} \, \exp[ \tilde{\lambda}_{\r} \, x ] \,
\beta_{\r} \right\} + \nonumber \\
 &  &
\sum_{\p \= \1}^{\n} \left\{
v_{\p} \, \cos \left[ b_{\p} \, x \right]
- \tilde{v}_{\p} \, \sin \left[ b_{\p} \, x \right] \right\} \,
\exp \left[  a_{\p} \, x  \right]  \, \gamma_{\p} + \nonumber \\
&  &
\sum_{\p \= \1}^{\n} \left\{
\tilde{v}_{\p} \, \cos \left[ b_{\p} \, x \right]
+ v_{\p} \, \sin \left[ b_{\p} \, x \right] \right\} \,
\exp \left[  a_{\p} \, x  \right]  \, \delta_{\p}~,
~~~\mbox{\small $m+n = 4$}~,
\end{eqnarray}
where $a_{\p} = \mbox{Re} \left[ z_{\p} \right]$,
 $b_{\p} = \mbox{Im} \left[ z_{\p} \right]$ and
$u_{\r}, \tilde{u}_{\r}, v_{\p}, \tilde{v}_{\p} \in \mathbb{H}$.
The $8$ real constants $\alpha_{\r}$, $\beta_{\r}$,
$\gamma_{\p}$, and $\delta_{\p}$  are fixed by the quaternionic initial
conditions $\varphi(x_{\0})=\varphi_{\0}$ and
$\frac{d \varphi}{dx} (x_{\0})= \varphi_{\1}$. The important point
to note here is that the particular solutions corresponding to
the complex eigenvalue $z=a+ib$ are given by a quaternionic
combination of  $\cos \left[ b \, x \right] \,
\exp \left[  a \, x  \right]$ and
$\sin \left[ b  \, x \right] \,
\exp \left[  a \, x  \right]$, namely
\[
\left\{ v \, \cos \left[ b \, x \right] -
\tilde{v} \, \sin \left[ b  \, x \right] \right\}
\exp \left[  a \, x  \right]~~~\mbox{and}~~~
\left\{ \tilde{v} \, \cos \left[ b \, x \right] +
v \, \sin \left[ b  \, x \right] \right\}
\exp \left[  a \, x  \right]~.
\]
For $\mathbb{C}$ linear quaternionic differential operators
(an explicit example is given in subsection~B),
the general solution
(\ref{genr}) reduces to
\begin{equation}
\varphi(x)  =  \sum_{\r \= \1}^{\m}
u_{\r} \, \exp[ \lambda_{\r} \, x ] \, c_{\r}  +
 \sum_{\p \= \1}^{\n}
v_{\p} \, \exp \left[  z_{\p} \, x  \right]  \, d_{\p}~,~~~
c_{\r} \, , \, d_{\p} \in \mathbb{C}~.
\end{equation}


\subsection*{A. Null nilpotent matrix}

Let us solve the second
order $\mathbb{R}$ linear quaternionic
differential  equation
\begin{equation}
\left[ \,
\frac{d^{^{\, \2}}}{dx^{^{\2}}} \, -
 L_{i} R_{j}
\, \frac{d}{dx} -
L_{j} R_{i} \, \right] \, \varphi(x) = 0~,
\end{equation}
with initial conditions $\varphi(0)=j $ and $\frac{d \varphi}{dx} (0)=k$.
The matrix
operator corresponding to this equation is
\begin{equation}
M_{\2}^{\Re} = \left( \begin{array}{cc} 0 & 1\\
                               L{\j} R_{\i} & ~L_{\i}  R_{\j}~
                 \end{array}
          \right)~.
\end{equation}
By using the matrix $T_{\2}^{\Re}$
\[
\begin{array}{lcl}
\{T^{\Re}_{\2}\}_{\1 \1}  & = &  \frac{1 + \sqrt{5}}{4} \,
\left( L_{\i} + L_{\k} R_{\j} -
L_{\i} R_{\k} - L_{\k} R_{\i}  \right)  -
\frac{1 - \sqrt{5}}{4} \, \left(
1  - L_{\j} R_{\j} + L_{\j} R_{\i} - R_{\k} \right)~,\\
\{T^{\Re}_{\2}\}_{\1 \2}  & = &
\frac{\sqrt{3}}{2 \sqrt{2}} \,
\left( 1 + L_{\i} R_{\i} - L_{\j} R_{\i} - L_{\k} \right) -
\frac{1}{2 \sqrt{2}} \,
\left( L_{\i} +  R_{\i}  + L_{\j} + L_{\k} R_{\i}  \right)~,\\
\{T^{\Re}_{\2}\}_{\2 \1}  & = & \frac{1}{2} \,
\left(  L_{\i} R_{\i}  + L_{\k} R_{\k} - L_{\j} R_{\k}
+ R_{\i}  - L_{\j}  + R_{\j} + L_{\k}   + L_{\i} R_{\j} \right)~,\\
\{T^{\Re}_{\2}\}_{\2 \2}  & = & \frac{1}{2 \sqrt{2}} \,
\left( L_{\j} R_{\k} - L_{\k} R_{\j}
- L_{\i} R_{\k} + R_{\j} \right)~,
\end{array}
\]
and its inverse
\[
\begin{array}{lcl}
\{(T^{\Re}_{\2})^{\mi \1}\}_{\1 \1} & = &
\frac{1}{4 \sqrt{5}} \, \left(
1  - L_{\j} R_{\j} + L_{\k} R_{\j}
- L_{\i} - L_{\k} R_{\i} - L_{\i} R_{\k} +
 L_{\j} R_{\i} + R_{\k} \right)~,\\
\{(T^{\Re}_{\2})^{\mi \1}\}_{\1 \2}  & = &
\frac{5 + \sqrt{5}}{40} \,
 \left( L_{\i} R_{\i}  + L_{\k} R_{\k} + L_{\i} R_{\j} - L_{\k}
\right) +  \frac{5 - \sqrt{5}}{40} \,
\left(
L_{\j}  - R_{\j} - L_{\j} R_{\k} - R_{\i} \right)~,\\
\{(T^{\Re}_{\2})^{\mi \1}\}_{\2 \1}  & = &
\frac{1}{2 \sqrt{6}} \,
\left(
1 + L_{\i} R_{\i} - L_{\j} R_{\i} + L_{\k}  \right)~,\\
\{(T^{\Re}_{\2})^{\mi \1}\}_{\2 \2}  & = &
\frac{1}{4 \sqrt{6}} \,
\left(
 L_{\j} R_{\j} - L_{\k} R_{\k} + L_{\i} R_{\j} + R_{\k} \right) +
\frac{1}{4 \sqrt{2}} \, \left(
 L_{\k} R_{\j} - L_{\j} R_{\k}  + L_{\i} R_{\k} + R_{\j} \right)~,
\end{array}
\]
we can rewrite $M_{\2}^{\Re}$ in terms of
\[
D_{\2}^{\Re} = \mbox{$\frac{1}{2}$} \, \left( \begin{array}{cc}
L_{\i} R_{\i} + \sqrt{5} \, L_{\k} R_{\k} & 0\\
                0  &  L_{\j}  R_{\j}~
                 \end{array}
          \right)~~~\mbox{and}~~~
A_{\2}^{\Re} = \mbox{$\frac{\sqrt{3}}{2}$} \,
            \left( \begin{array}{cc} 0 & 0\\
                           0  &   R_{\j}
                 \end{array}
          \right)~.
\]
From Eq.~(\ref{matsol}) we obtain the following solution
\begin{eqnarray}
\varphi (x) & = & \{T^{\Re}_{\2}\}_{\1 \1} \,
\exp \left[ \, \mbox{$\frac{L_{\i} R_{\i} + \sqrt{5} \, L_{\k} R_{\k}}{2}$}
\, x \, \right] \, \left[
\{(T^{\Re}_{\2})^{\mi \1}\}_{\1 \1} \, \, j +
\{(T^{\Re}_{\2})^{\mi \1}\}_{\1 \2} \, \, k \right] +
\nonumber \\
 & &
\{T^{\Re}_{\2}\}_{\1 \2} \, \exp \left[ \,
\mbox{$\frac{L_{\j} R_{\j} + \sqrt{3} \, R_{\j}}{2}$}
\, x \, \right]
 \, \left[
\{(T^{\Re}_{\2})^{\mi \1}\}_{\2 \1}\, \, j +
\{(T^{\Re}_{\2})^{\mi \1}\}_{\2 \2} \, \, k \right]
\nonumber \\
 & = &
\left\{ \mbox{$\frac{i - j}{2 \sqrt{5}}$} \, \left(
\sinh \left[\mbox{$\frac{\sqrt{5}}{2}$} \, x \right] - \sqrt{5} \,
\cosh \left[\mbox{$\frac{\sqrt{5}}{2}$} \, x \right] \right)
+ \mbox{${\frac{1 + k}{\sqrt{3}}}$} \,
\sin \left[\mbox{$\frac{\sqrt{3}}{2}$} \, x \right] \right\} \,
\exp \left[ x/2  \right] +   \nonumber\\
 & &
\left\{ \mbox{$\frac{k - 1}{\sqrt{5}}$} \,
\sinh \left[\mbox{$\frac{\sqrt{5}}{2}$} \, x \right]
+ \mbox{$\frac{i + j}{2}$} \, \left(
\cos \left[\mbox{$\frac{\sqrt{3}}{2}$} \, x \right]
+ \mbox{${\frac{1}{\sqrt{3}}}$} \,
\sin \left[\mbox{$\frac{\sqrt{3}}{2}$} \right] \right)
\right\} \,
\exp \left[ - x/2  \right]~.
\end{eqnarray}
In the case of  a null nilpotent matrix,
the general solution of real linear quaternionic differential equations
can be written in terms of  the eigenvalue spectrum of
$M_{\8}[\mathbb{R}]$, real matrix counterpart  of $M_{\2}^{\Re}$.
In this particular case, the eigenvalue spectrum is
\[
\left\{
\mbox{$\frac{1+\sqrt{5}}{2}$} \, , \,
\mbox{$\frac{1 - \sqrt{5}}{2}$} \, , \,
- \mbox{$\frac{1+\sqrt{5}}{2}$} \, , \,
- \mbox{$\frac{1 - \sqrt{5}}{2}$} \, , \,
\mbox{$\frac{1+ i \sqrt{3}}{2}$} \, , \,
- \mbox{$\frac{1 - i\sqrt{3}}{2}$} \, , \,
\mbox{$\frac{1 - i \sqrt{3}}{2}$} \, , \,
- \mbox{$\frac{1 + i \sqrt{3}}{2}$} \, \right\}~.
\]
Consequently, Eq.~(\ref{genr}) becomes
\begin{eqnarray*}
\varphi(x) =
u_{\1} \,
\exp \left[ \mbox{$\frac{1+\sqrt{5}}{2}$} \, x \right] \, \alpha_{\1} +
u_{\2} \,
\exp \left[ \mbox{$\frac{1 - \sqrt{5}}{2}$} \, x \right] \, \alpha_{\2}+
 \tilde{u}_{\1} \,
\exp \left[ - \mbox{$\frac{1+\sqrt{5}}{2}$} \, x \right] \, \beta_{\1} +
\tilde{u}_{\2} \,
\exp \left[ - \mbox{$\frac{1 - \sqrt{5}}{2}$} \, x \right] \, \beta_{\2}
 + \\
\left\{ v_{\1} \, \cos \left[ \mbox{$\frac{\sqrt{3}}{2}$} \, x \right]
- \tilde{v}_{\1} \, \sin \left[ \mbox{$\frac{\sqrt{3}}{2}$}
\, x \right] \right\} \,
\exp \left[   x /2   \right]  \, \gamma_{\1} +
\left\{ \tilde{v}_{\1} \, \cos \left[ \mbox{$\frac{\sqrt{3}}{2}$} \, x \right]
+ v_{\1} \, \sin \left[ \mbox{$\frac{\sqrt{3}}{2}$}
\, x \right] \right\} \,
\exp \left[   x /2   \right]  \, \delta_{\1} + \\
\left\{ v_{\2} \, \cos \left[ \mbox{$\frac{\sqrt{3}}{2}$} \, x \right]
- \tilde{v}_{\2} \, \sin \left[ \mbox{$\frac{\sqrt{3}}{2}$}
\, x \right] \right\} \,
\exp \left[   - x /2   \right]  \, \gamma_{\2} +
\left\{ \tilde{v}_{\2} \, \cos \left[ \mbox{$\frac{\sqrt{3}}{2}$} \, x \right]
+ v_{\2} \, \sin \left[ \mbox{$\frac{\sqrt{3}}{2}$}
\, x \right] \right\} \,
\exp \left[   - x /2   \right]  \, \delta_{\2}~.~
\end{eqnarray*}
By direct calculations, we can determine the quaternionic coefficients
$u_{\1,\2}$, $\tilde{u}_{\1,\2}$, $v_{\1,\2}$ and $\tilde{v}_{\1,\2}$. We find
\[
\begin{array}{lcl}
u_{\1} = i - j &~,~~~& \tilde{u}_{\1}=1-k~,\\
u_{\2} = i - j &~,~~~& \tilde{u}_{\2}=1-k~,\\
v_{\1} = 1 + k  &~,~~~& \tilde{v}_{\1}=0~,\\
v_{\2} = i + j &~,~~~& \tilde{v}_{\2}=0~.
\end{array}
\]
The particular solutions corresponding to the complex eigenvalues of
$M_{\8}[\mathbb{R}]$ are thus given by
\[
z_{\1} = \mbox{$\frac{1+ i \sqrt{3}}{2}$} \, , \, \bar{z}_{\1}~:~~~
(1 + k) \, \exp \left[  x /2   \right] \, ~\begin{array}{c}
\cos \left[ \mbox{$\frac{\sqrt{3}}{2}$} \, x \right] \\~\\
 \sin
\left[ \mbox{$\frac{\sqrt{3}}{2}$} \, x \right]
\end{array}~,
\]
and
\[
z_{\2} = - \, \mbox{$\frac{1- i \sqrt{3}}{2}$} \, , \, \bar{z}_{\2}
~:~~~
(i+j) \, \exp \left[ - \,   x /2   \right] \, ~\begin{array}{c}
\cos \left[ \mbox{$\frac{\sqrt{3}}{2}$}\right] \\~\\
 \sin
\left[ \mbox{$\frac{\sqrt{3}}{2}$}\right]
\end{array}~,
\]
The  initial conditions
$\varphi(0)=j$ and $\frac{d \varphi}{dx} (0)=k$ make the solution completely
determined,
\[
\begin{array}{lcrclcrclccclccr}
\alpha_{\1} & = & \frac{1- \sqrt{5}}{4 \sqrt{5}} &~,~~&
\beta_{\1} & =  & \frac{1}{2\sqrt{5}}&~,~~&
\gamma_{\1} &= & 0 &~,~~&
\delta_{\1}  & = &  \frac{1}{\sqrt{3}} &~,\\
\alpha_{\2} & = & -\frac{1+  \sqrt{5}}{4 \sqrt{5}} &~,~~&
\beta_{\2} & = & -\frac{1}{2\sqrt{5}} &~,~~&
\gamma_{\2} & = & \frac{1}{2}  &~,~~&
\delta_{\2} & = & \frac{1}{2 \sqrt{3}} &~.
\end{array}
\]

Let us now consider the following second
order $\mathbb{R}$ linear quaternionic
differential  equation
\begin{equation}
\label{eqv}
\left[ \,
\frac{d^{^{\, \2}}}{dx^{^{\2}}} \, -
\left(  L_{i} R_{i} + L_{\j} R_{\j} \right)
\, \frac{d}{dx} +
R_{k} -  L_{k} \, \right] \, \varphi(x) = 0~.
\end{equation}
The eigenvalue spectrum of the  real matrix counterpart  of
\begin{equation}
M_{\2}^{\Re} = \left( \begin{array}{cc} 0 & 1\\
           ~L{\k} - R_{\k}~ & ~L_{\i}  R_{\i} + L_{\j} R_{\j}~
                 \end{array}
          \right)
\end{equation}
is
\[
\left\{
2 \, , \, - 2 \, , \,
0 \, , \,
0 \, , \,
1 + i
\, , \,
- (1 - i) \, , \,
1 - i  \, , \,
- (1 + i)  \, \right\}~.
\]
In this case $v_{\1} \cos x \, \exp[x]$
and $v_{\2} \cos x \, \exp[-x]$ does not represent (non trivial)
particular solutions.
In fact, in order to satisfy Eq.~(\ref{eqv}) we have to impose the following
constraints on the quaternionic coefficients $v_{\1,\2}$
\[
\begin{array}{rcl}
\left[ L_{\i}R_{\i} + L_{\j}R_{\j} - 2 \right] \, v_{\1} = 0 &
~~\Rightarrow~~ & v_{1} = i \, \beta + j \, \gamma~,\\
\left[ L_{\i}R_{\i} + L_{\j}R_{\j} + L_{\k} - R_{\k} \right] \, v_{\1} = 0 &
~~\Rightarrow~~ & v_{1} = 0~,\\
\left[ L_{\i}R_{\i} + L_{\j}R_{\j} + 2 \right] \, v_{\2} = 0 &
~~\Rightarrow~~ & v_{2} = \alpha + j \, \delta~,\\
\left[ L_{\i}R_{\i} + L_{\j}R_{\j} - L_{\k} + R_{\k} \right] \, v_{\2} = 0 &
~~\Rightarrow~~ & v_{2} = 0~.
\end{array}
\]
The particular solutions corresponding to the complex eigenvalues
$z_{\1} = 1+i$ and
$z_{\2} = -(1 +i)$ are
\[
z_{\1} \, , \, \bar{z}_{\1}~:~~~
\exp \left[  x   \right] \, ~\begin{array}{c}
 j \, \cos \left[ x \right] - i \, \sin \left[ x \right] \\~\\
 i \, \cos \left[ x \right] + j \, \sin \left[ x \right]
\end{array}~,
\]
and
\[
z_{\2} \, , \, \bar{z}_{\2}
~:~~~
\exp \left[  -  x   \right] \, ~\begin{array}{c}
 j \, \cos \left[ x \right] + i \, \sin \left[ x \right] \\~\\
 i \, \cos \left[ x \right] - j \, \sin \left[ x \right]
\end{array}~.
\]
The general solution is given by
\begin{eqnarray*}
\varphi(x) & = &
\exp \left[ - \, 2 \, x \right] \, \alpha_{\1} +
k \,
\exp \left[ 2 \, x \right] \, \alpha_{\2} +
\beta_{\1} +
k \, \beta_{\2}
 + \\
& &
\left\{ j \, \cos \left[ x \right]
- i \, \sin \left[ x \right] \right\} \,
\exp \left[   x    \right]  \, \gamma_{\1} +
\left\{ i \, \cos \left[ x \right]
+ j \, \sin \left[ x \right] \right\} \,
\exp \left[   x   \right]  \, \delta_{\1} + \\
& &
\left\{ j \, \cos \left[ x \right]
+ i \, \sin \left[ x \right] \right\} \,
\exp \left[   -  x    \right]  \, \gamma_{\2} +
\left\{ i \, \cos \left[ x \right]
-  j \, \sin \left[ x \right] \right\} \,
\exp \left[   - x    \right]  \, \delta_{\2}~.~
\end{eqnarray*}

\subsection*{B. Complex linear case}

We now determine the solution of the second
order $\mathbb{C}$ linear quaternionic
differential  equation
\begin{equation}
\left[ \,
\frac{d^{^{\, \2}}}{dx^{^{\2}}} \, -
L_{j} R_{\i} \, \right] \, \varphi(x) = 0~.
\end{equation}
The eigenvalue spectrum of the  real matrix counterpart  of
\begin{equation}
M_{\2}^{\Re} = \left( \begin{array}{cc} 0 & 1\\
           ~L{\j} R_{\i}~ & 0
                 \end{array}
          \right)
\end{equation}
is
\[
\left\{
1 \, , \, -1 \, , \,
1 \, , \,
-1 \, , \,
i
\, , \,
i \, , \,
- i  \, , \,
- i  \, \right\}~.
\]
The general solution is then given by
\begin{eqnarray}
\label{eqc}
\varphi(x) & = &
u_{\1} \,
\exp \left[ x \right] \, \alpha_{\1} +
u_{\2} \,
\exp \left[ -  x \right] \, \alpha_{\2}+
 \tilde{u}_{\1} \,
\exp \left[  x \right] \, \beta_{\1} +
\tilde{u}_{\2} \,
\exp \left[ -  x \right] \, \beta_{\2}
 + \nonumber \\
 & &
\left\{ v_{\1} \, \cos \left[ x \right] -
\tilde{v}_{\1} \, \sin \left[  x \right] \right\} \, \gamma_{\1} +
\left\{ \tilde{v}_{\1} \, \cos \left[ x \right]
+ v_{\1} \, \sin \left[ x \right] \right\} \, \delta_{\1} +
\nonumber \\
& & \left\{ v_{\2} \, \cos \left[ x \right]
- \tilde{v}_{\2} \, \sin \left[ x \right] \right\} \, \gamma_{\2} +
\left\{ \tilde{v}_{\2} \, \cos \left[ x \right]
+ v_{\2} \, \sin \left[x \right] \right\} \, \delta_{\2}~.~
\end{eqnarray}
By direct calculations, we can determine the quaternionic coefficients
$u_{\1,\2}$, $\tilde{u}_{\1,\2}$, $v_{\1,\2}$ and $\tilde{v}_{\1,\2}$.
We find
\[
\begin{array}{lcl}
u_{\1} = k - 1 &~,~~~& \tilde{u}_{\1}=j-i~,\\
u_{\2} = k - 1 &~,~~~& \tilde{u}_{\2}=j-i~,\\
v_{\1} = 1 + k  &~,~~~& \tilde{v}_{\1}=0~,\\
v_{\2} = i + j &~,~~~& \tilde{v}_{\2}=0~.
\end{array}
\]
Consequently, Eq.~(\ref{eqc}) becomes
\begin{equation}
\varphi(x) = \left(k - 1 \right) \, \left[
\exp \left[ x \right] \, c_{\1} +
\exp \left[ - x \right] \, c_{\2} \right] +
\left(1 + k  \right) \, \left\{
\exp \left[ i  x \right] \, d_{\1}+
\exp \left[ - i   x \right] \, d_{\2}  \right\}~.
\end{equation}
The initial conditions $\varphi(0)=j $ and $\frac{d \varphi}{dx} (0)=k$
fix the following solution
\begin{equation}
\varphi(x) = \left\{ (j - i)\,
\cosh x  + (k - 1) \,
\sinh x  + (i + j) \,
\exp \left[ - i x \right] \right\} / 2~.
\end{equation}

\subsection*{C. Non null nilpotent matrix}

As a last example we consider the  differential
equation
\begin{equation}
\left[ \,
\frac{d^{^{\, \2}}}{dx^{^{\2}}} \, -
\left(  L_{i} R_{j} + L_{\j} R_{\i} \right)
 \right] \, \varphi(x) = 0~.
\end{equation}
The eigenvalue spectrum of the  real matrix counterpart  of
\begin{equation}
M_{\2}^{\Re} = \left( \begin{array}{cc} 0 & 1\\
           ~L_{\i}  R_{\j} + L_{\j} R_{\i}~ & 0
                 \end{array}
          \right)
\end{equation}
is
\[
\left\{
0 \, , \, 0 \, , \, 0 \, , \, 0 \, , \, \sqrt{2} \, , \,
- \sqrt{2} \, , \, i \sqrt{2} \, , \, - i \sqrt{2}
\, \right\}~.
\]
By using an appropriate matrix $T_{\2}^{\Re}$, we can rewrite
$M_{\2}^{\Re}$ in terms of
\[
D_{\2}^{\Re} = \mbox{$\frac{1}{\sqrt{2}}$} \, \left( \begin{array}{cc}
 0  & 0\\
                0  &  L_{\j} R_{\j} + L_{\k} R_{\k}
                 \end{array}
          \right)~,~~
A_{\2}^{\Re} = \mbox{$\frac{1}{\sqrt{2}}$} \,
            \left( \begin{array}{cc} 0 & 0\\
                           0  & L_{\i} -   R_{\i}
                 \end{array}
          \right)~~\mbox{and}~~
N_{\2}^{\Re} = \mbox{$\frac{1}{2}$} \,
            \left( \begin{array}{cc}  L_{\i} - L_{\k} R_{\j} & 0\\
                           0  &   0
                 \end{array}
          \right)~,
\]
$M_{\2}^{\Re} = T_{\2}^{\Re} \left( D_{\2}^{\Re} +
A_{\2}^{\Re} + N_{\2}^{\Re}  \right)
(T_{\2}^{\Re})^{\mi \1}$. Notwithstanding its diagonal form,
$N_{\2}^{\Re}$ represents a nilpotent matrix. In fact, an easy algebraic
calculation shows that $(N_{\2}^{\Re})^{\2}=0$. This implies in the general
solution the explicit presence of the real variable $x$. Explicitly, we
have
\begin{eqnarray*}
\varphi(x) & = &
\alpha_{\1} +
k \, \alpha_{\2} +
\left( \tilde{\alpha}_{\1} +
k \, \tilde{\alpha}_{\2} \right) \, x + \\
 & & (i - j) \, \left\{
\exp \left[  \sqrt{2} \,  x   \right] \, \beta_{\1} +
\exp \left[ -  \sqrt{2} \,  x   \right] \, \beta_{\2} \right\}
 + \\
& &
(i+j) \,
\left\{ \cos \left[ \sqrt{2} \, x \right] \gamma_{\1} +
         \sin \left[ \sqrt{2} \, x \right] \, \delta_{\1} \right\}~.
\end{eqnarray*}
To make this solution determined, let us specify certain supplementary
constraints, $\varphi(0)=j$ and $\frac{d \varphi}{dx} (0)=k$,
\begin{equation}
\varphi (x) = k x + \mbox{$\frac{j - i }{2}$} \,
\cosh \left[ \sqrt{2} \, x \right] +
\mbox{${\frac{i + j}{2}}$} \,
\cos \left[ \sqrt{2} \, x \right]~.
\end{equation}
By changing the initial conditions,
$\varphi(0) = k$ and $\frac{d \varphi}{dx} (0)=j$, we find
\begin{equation}
\varphi(x) =  k  + \mbox{$\frac{j - i }{2\sqrt{2}}$} \,
\sinh \left[ \sqrt{2} \, x \right] +
\mbox{${\frac{i + j}{2 \sqrt{2}}}$} \,
\sin \left[ \sqrt{2} \, x \right]~.
\end{equation}


\section*{III. Conclusions}

In this paper, we have extended the resolution method of $\mathbb{H}$
and $\mathbb{C}$ linear ordinary differential equations with quaternionic
constant coefficients~\cite{DeDu01} to $\mathbb{R}$ linear differential
operators. By a matrix approach, we have shown that particular solutions
of differential equations with constant quaternionic coefficients which
appear both on the left and on the right can be given in terms of the
eigenvalues of the matrix
$M_{\4 \n}[\mathbb{R}]$, representing the real counterpart of
the quaternionic operators, $M_{\n}^{\Re}$, associated to the $\mathbb{R}$
linear quaternionic differential equation. In correspondence to the
eigenvalue $z=a+ib$, we have found the following particular solutions
\begin{equation}
\label{parsol}
\left\{ q \, \cos \left[ b \, x \right] -
p \, \sin \left[ b  \, x \right] \right\}
\exp \left[  a \, x  \right]~~~\mbox{and}~~~
\left\{ p \, \cos \left[ b \, x \right] +
q \, \sin \left[ b  \, x \right] \right\}
\exp \left[  a \, x  \right]~,~~~q \, , \, p \in \mathbb{H}~.
\end{equation}
When $q=p$ or $p=0$ these solutions reduce to
\[
q \, \cos \left[ b \, x \right] \,
\exp \left[  a \, x  \right]~~~\mbox{and}~~~
q  \, \sin \left[ b  \, x \right] \,
\exp \left[  a \, x  \right]~.
\]
For $\mathbb{C}$ linear differential operators, the particular
solutions~(\ref{parsol}) couple to give
\[
q \, \exp \left[  (a  + i b) \, x  \right]~.
\]
Our discussion can be viewed as a preliminary step towards a full
understanding of the role that quaternions could play in analysis,
linear algebra and, consequently, in physical applications.
A complete theory of quaternionic differential operators (as well of the
quaternionic eigenvalue problem) is at present far from being conclusive
and deserves further investigations.

\section*{Acknowledgements}

S.~D.~L. wishes to thank Prof.~L.~Solombrino and Dr.~G.~Scolarici of the
Department of Physics of the University of Lecce for comments and
discussions about the eigenvalue problem for real linear quaternionic
operators. G.~C.~D. is grateful to CAPES for financial support.



\newpage


\section*{Appendix. Real translation tables}

\vspace*{0.5cm}

\noindent {\bf Table 1}.

{\footnotesize
\[
\hspace*{-1cm}
\begin{array}{c||ccc}
 ~~\left(
    \begin{array}{cccc}
     ~1~ & ~0~ & ~0~ & ~0~\\
     ~0~ & 1 & 0 & 0\\
     0 & 0 & 1 & 0\\
     0 & 0 & 0 & 1
   \end{array}  \right)~~   &  ~~\left(
     \begin{array}{cccc}
       0 & $-$1\, \, & 0 & 0\\
       ~1~ & ~0~ & ~0~ & ~0~\\
       ~0~ & 0 & 0 & 1\\
       0 & 0 & $-$1\, \, & 0
     \end{array}  \right)~~   &  \left(
        \begin{array}{cccc}
          0 & 0 & $-$1\, \, & 0\\
          ~0~ & 0 & 0 & $-$1\, \,\\
          ~1~ & ~0~ & ~0~ & ~0~\\
          0 & 1 & 0 & 0
        \end{array}  \right)   &  ~~\left(
           \begin{array}{cccc}
             0 & 0 & 0 & $-$1\, \,\\
             ~0~ & ~0~ & ~1~ & ~0~\\
             0 & $-$1\, \, & ~0~ & 0\\
             1 & 0 & 0 & 0
           \end{array}  \right)~~ \\
1 & R_{\i} & R_{\j} & R_{\k} \\ [1mm]
\hline \hline
 & & & \\
  \left(
    \begin{array}{cccc}
     0 & $-$1\, \, & 0 & 0\\
     ~1~ & ~0~ & ~0~ & ~0~\\
     ~0~ & 0 & 0 & $-$1\, \,\\
     0 & 0 & 1 & 0
   \end{array}  \right)  &  \left(
     \begin{array}{cccc}
       $-$1\, \, & 0 & 0 & 0\\
       0 & $-$1\, \, & 0 & 0\\
       ~0~ & ~0~ & ~1~ & ~0~\\
       0 & 0 & ~0~ & 1
     \end{array}  \right)   &  \left(
        \begin{array}{cccc}
          ~0~ & ~0~ & ~0~ & ~1~\\
          0 & 0 & $-$1\, \, & ~0~ \\
          0 & $-$1\, \, & 0 & 0\\
          1 & 0 & 0 & 0
        \end{array}  \right)  &  \left(
           \begin{array}{cccc}
             0 & 0 & $-$1\, \, & 0\\
             0 & 0 & 0 & $-$1\, \,\\
             $-$1\, \, & ~0~ & ~0~ & ~0~\\
             ~0~ & $-$1\, \, & 0 & 0
           \end{array}  \right) \\
L_{\i} & L_{\i}\,R_{\i} & L_{\i}\,R_{\j} & L_{\i}\,R_{\k} \\
 & & & \\
   \left(
    \begin{array}{cccc}
     0 & 0 & $-$1\, \, & 0\\
     ~0~ & ~0~ & ~0~ & ~1~\\
     1 & 0 & 0 & ~0~\\
     0 & $-$1\, \, & 0 & 0
   \end{array}  \right)  &  \left(
     \begin{array}{cccc}
       0 & 0 & 0 & $-$1\, \,\\
       0 & 0 & $-$1\, \, & 0\\
       ~0~ & $-$1\, \, & 0 & 0\\
       $-$1\, \, & ~0~ & ~0~ & ~0~
     \end{array}  \right)   &  \left(
        \begin{array}{cccc}
          $-$1\, \, & 0 & 0 & 0\\
          ~0~ & ~1~ & ~0~ & ~0~\\
          0 & ~0~ & $-$1\, \, & 0\\
          0 & 0 & 0 & 1
        \end{array}  \right)  &  \left(
           \begin{array}{cccc}
             ~0~ & ~1~ & ~0~ & ~0~\\
             1 & ~0~ & 0 & 0\\
             0 & 0 & 0 & $-$1\, \,\\
             0 & 0 & $-$1\, \, & 0
           \end{array}  \right) \\
L_{\j} & L_{\j}\,R_{\i} & L_{\j}\,R_{\j} & L_{\j}\,R_{\k} \\
  & & & \\
   \left(
    \begin{array}{cccc}
     0 & 0 & 0 & $-$1\, \,\\
     0 & 0 & $-$1\, \, & 0\\
     ~0~ & ~1~ & ~0~ & ~0~\\
     1 & ~0~ & 0 & 0
   \end{array}  \right) & \left(
     \begin{array}{cccc}
       ~0~ & ~0~ & ~1~ & ~0~\\
       0 & 0 & 0 & $-$1\, \,\\
       1 & 0 & ~0~ & 0\\
       0 & $-$1\, \, & 0 & 0
     \end{array}  \right)   &  \left(
        \begin{array}{cccc}
          0 & $-$1\, \, & 0 & 0\\
          $-$1\, \, & 0 & 0 & 0\\
          ~0~ & ~0~ & ~0~ & $-$1\, \,\\
          0 & 0 & $-$1\, \, & ~0~
        \end{array}  \right)   &  \left(
           \begin{array}{cccc}
             $-$1\, \, & 0 & 0 & 0\\
             ~0~ & ~1~ & ~0~ & ~0~\\
             0 & ~0~ & 1 & 0\\
             0 & 0 & 0 & $-$1\, \,
           \end{array}  \right) \\
L_{\k} & L_{\k}\,R_{\i} & L_{\k}\,R_{\j} & L_{\k}\,R_{\k}
\end{array}
\]
}

\vspace*{1cm}

\noindent {\bf Example}.

\begin{eqnarray*} \hspace*{-1cm}
L_{\i} -
2 \, R_{\j} - 3 \, L_{i} R_{\k} & ~\rightarrow~ &
 \left(
           \begin{array}{cccc}
             0 & $-$1  & 0 & 0\\
            ~1~ & ~0~ & ~0~ & ~0~\\
             ~0~ & 0 & 0 & $-$1\, \,\\
             0 & 0 & 1 & 0
           \end{array}  \right)
- 2 \,
 \left(
           \begin{array}{cccc}
             0 & 0 & $-$1\, \, & 0\\
             ~0~ & ~0~ & ~0~ & $-$1\, \,\\
             1 & 0 & 0 & ~0~\\
             0 & 1 & 0 & 0
           \end{array}  \right)
- 3 \,
\left(
           \begin{array}{cccc}
             0 & 0 & $-$1\, \, & 0\\
             ~0~ & ~0~ & ~0~ & $-$1\, \,\\
             $-$1\, \, & 0 & 0 & ~0~\\
             0 & $-$1\, \, & 0 & 0
           \end{array}  \right)\\
 & ~\rightarrow~ &
\left(
           \begin{array}{cccc}
             0  & $-$1\, \, & 5 & 0\\
             1 & ~0~ & ~0~ & 5\\
             1 & ~0~ & ~0~ & $-$1\, \,\\
             ~0~ & 1 & 1 & ~0~
           \end{array}  \right)~.
\end{eqnarray*}

\newpage

\noindent {\bf Table 2}.

\[
\begin{array}{ccccc}
\begin{array}{||c||c|c|c|c||}
\hline\hline
~1/4~  & ~~1~~ & L_{\i}R_{\i} & L_{\j}R_{\j} & L_{\k}R_{\k}\\
\hline\hline
a_{\0\0} &  + &  -  &  -  &  -  \\
a_{\1\1} &  + &  -  &  +  &  +  \\
a_{\2\2} &  + &  +  &  -  &  +  \\
a_{\3\3} &  + &  +  &  +  &  -  \\
\hline\hline
 ~1/4~    &  ~L_{\i}~ & ~R_{\i}~ & L_{\j}R_{\k} & L_{\k}R_{\j}\\
\hline \hline
a_{\0\1} &  - &  -  &  +  &  - \\
a_{\1\0} &  + &  +  &  +  &  - \\
a_{\2\3} &  - &  +  &  -  &  -  \\
a_{\3\2} &  + &  -  &  -  &  -   \\
\hline\hline
\end{array}  &  & &
\begin{array}{||c||c|c|c|c||}
\hline\hline
~1/4~   & ~L_{\j}~ & ~R_{\j}~ & L_{\i}R_{\k} & L_{\k}R_{\i}\\
\hline \hline
a_{\0\2} & - &  -  &  -  &  + \\
a_{\2\0} & + &  +  &  -  &  + \\
a_{\1\3} & + &  -  &  -  &  - \\
a_{\3\1} & - &  +  &  -  &  -  \\
\hline\hline
~1/4~  & ~L_{\k}~ & ~R_{\k}~ & L_{\i}R_{\j} & L_{\j}R_{\i}\\
\hline \hline
a_{\0\3} & - &  -  &  +  &  - \\
a_{\3\0} & + &  +  &  +  &  - \\
a_{\1\2} & - &  +  &  -  &  - \\
a_{\2\1} & + &  -  &  -  &  - \\
\hline\hline
\end{array}
\end{array}
\]

\vspace*{1cm}

\noindent {\bf Example}.
\[
A :=
\left(
           \begin{array}{cccc}
             ~a_{\0 \0}~  & ~a_{\0 \1}~ & ~a_{\0 \2}~ & ~a_{\0 \3}~ \\
             a_{\1 \0}  & a_{\1 \1} & a_{\1 \2} & a_{\1 \3} \\
             a_{\2 \0}  & a_{\2 \1} & a_{\2 \2} & a_{\2 \3} \\
             a_{\3 \0}  & a_{\3 \1} & a_{\3 \2} & a_{\3 \3}
             \end{array}  \right) =
\left(
           \begin{array}{cccc}
             0  & $-$1\, \, & 5 & 0\\
             1 & ~0~ & ~0~ & 5\\
             1 & ~0~ & ~0~ & $-$1\, \,\\
             ~0~ & 1 & 1 & ~0~
           \end{array}  \right)~.\]
 \begin{eqnarray*}
A & ~\rightarrow~&
( a_{\0 \1} / 4 ) \, \left[
- L_{\i} -  R_{\i} +  L_{\j}R_{\k} - L_{\k}R_{\j}
 \right]+\\ & &
( a_{\0 \2} / 4 ) \, \left[
- L_{\j} -  R_{\j} -  L_{\i}R_{\k} + L_{\k}R_{\i}
\right]+\\ & &
( a_{\1 \0} / 4 ) \, \left[
+ L_{\i} +  R_{\i} +  L_{\j}R_{\k} - L_{\k}R_{\j}
\right]+\\ & &
( a_{\1 \3} / 4 ) \, \left[
+ L_{\j} -  R_{\j} -  L_{\i}R_{\k} - L_{\k}R_{\i}
\right]+\\ & &
( a_{\2 \0} / 4 ) \, \left[
+ L_{\j} +  R_{\j} -  L_{\i}R_{\k} + L_{\k}R_{\i}
\right]+\\ & &
( a_{\2 \3} / 4 ) \, \left[
- L_{\i} +  R_{\i} -  L_{\j}R_{\k} - L_{\k}R_{\j}
\right]+\\ & &
( a_{\3 \1} / 4 ) \, \left[
- L_{\j} +  R_{\j} -  L_{\i}R_{\k} - L_{\k}R_{\i}
\right]+\\ & &
( a_{\3 \2} / 4 ) \, \left[
+ L_{\i} -  R_{\i} +  L_{\j}R_{\k} - L_{\k}R_{\j}
\right]\\
& ~\rightarrow~ &
L_{\i} - 2 \, R_{\j} - 3 \, L_{i} R_{\k}~.
\end{eqnarray*}


\begin{thebibliography}{99}

\bibitem{DeDu01}
De Leo S and Ducati G C,
``Quaternionic differential operators'',
J.~Math.~Phys. {\bf 42}, 2236--2265 (2001).

\bibitem{Zha97}
Zhang F,
``Quaternions and matrices of quaternions'',
Lin.~Alg.~Appl. {\bf 251} 21--57 (1997).


\bibitem{Bak99}
Baker A,
``Right eigenvalues for quaternionic matrices: a topological approach'',
Lin.~Alg.~Appl. {\bf 286} 303--309 (1999).

\bibitem{DeSco00}
De Leo S and Scolarici G,
``Right eigenvalue equation in quaternionic quantum mechanics'',
J.~Phys.~A {\bf 33} 2971--2995 (2000).



 \end{thebibliography}
\end{document}